\documentclass{ws-procs9x6}
\usepackage{graphicx}
\usepackage{color}

\newcommand{\lt}{\triangleleft}

\newcommand{\Z}{\mathbb{Z}}
\newcommand{\R}{\mathbb{R}}

\newcommand{\Aut}{\mbox{Aut}}

\newcommand{\Inn}{\mbox{Inn}}

\usepackage{graphicx}

\def\Aut{{\rm Aut}}
\begin{document}

\title{A Survey of Quandle Ideas}
\author{
J. Scott Carter\footnote{Supported in part by
NSF Grant DMS  \#0603926.}
\\ University of South Alabama \\ carter@jaguar1.usouthal.edu
}

\begin{abstract} This article surveys many aspects of the theory of quandles which algebraically encode the Reidemeister moves. In addition to knot theory, quandles have found applications in other areas which are only mentioned in passing here. The main purpose is to give a short introduction to the subject and a guide to the applications that have been found thus far for quandle cocycle invariants. 
\end{abstract}

\section{Introduction}
A quandle is an algebraic system whose axioms are derived from the Reidemeister moves to classical knot diagrams. Any codimension-2, piecewise-linear, embedding whose image in Euclidean space is locally-flat has a fundamental quandle. In particular, classical knots are completely classified by their fundamental quandles. Of course, this means that the fundamental quandle is extremely hard to compute. In fact, I only know of one example (that of the trefoil) in which the fundamental quandle has been identified to a heretofore unrelated quandle (the Dehn quandle of simple closed curves on a torus in which the quandle operation is induced by Dehn twists\cite{NP4}).  The fundamental quandle  of a classical knot describes  the fundamental group and a peripheral subgroup, and it is in this respect that it classifies the knot. Section~\ref{Joyofsex} contains this description.

Interesting and computable invariants come from representing the fundamental quandle into a finite quandle that is often, but not always, a dihedral quandle. Not only can we count the number of such quandle colorings, but we can separate different colorings into distinct subclasses that are indexed by the value of a quandle cocycle. The quandle cocycle invariants are described below. 

The final section of the paper summarizes many of the applications that have been found for quandle cocycles. The list is certainly not encompassing, but contains, to the best of my   knowledge, a reasonable guide to the myriad results that have been obtained in the wake of the definition the quandle cocycle invariants\cite{CJKLS}. 

First I turn to a synopsis of the history of the subject and some important definitions. My historical perspective will be tainted by my vision based upon contributions to the field that my collaborators and I have made.

\subsection*{History, Definitions, and Examples}

 Joyce's dissertation cites the 1929 work of Burstin and Mayer \cite{BurMay} as the earliest work on distributive groups. He points out that similar ideas underlie the core group of Moufang loops \cite{BB,Moufang}.  
In 1941/2, Takasaki \cite{Taka} introduced the notion of a {\it kei} --- a set that possesses an involutory self-distributive operation. A 1959 correspondence between Gavin Wraith and John Conway re-introduced the idea in the context of a group that acts upon itself via conjugation. Therein the idea was named a wrack (after wrack and ruin of a group after the multiplication operation had been dismissed). According to Wikipedia and other sources, the name ``wrack" was a bit of Conway word play upon Wraith's name. Modern authors have returned to the spelling ``rack."  

In the terminology that was later introduced by Joyce \cite{Joyce}, a kei is an involutory quandle. At approximately the same time of Joyce, Matveev \cite{Matveev} introduced self-distributive groupoids. Matveev and Joyce showed that these structures determine the complement of a knot up to an orientation reversing homeomorphism. The date 1979 is common to both Joyce's thesis and Matveev's original talks. These works, it seems, were developed independently and at roughly the same time.  
 Steve Winker's dissertation was completed in 1982 under the direction of Lou Kauffman. In the dissertation, Winker shows that the fundamental involutory quandle is a complete classifier.

Brieskorn's vast article \cite{Bries} appeared in the proceedings of a 1986 Santa Cruz conference on Braids \cite{braids}. It is notable for a couple of reasons. First, the volume itself marks the rebirth of braids and the braid group as important tools at the interface among algebra, topology, and analysis. This interface was due, in part, because of Jones's seminal discovery \cite{Jones}. Second, Brieskorn's paper lays important foundations in which automorphic sets also play important roles among various branches of mathematics. Many of these connections remain unexplored. I encourage young researchers to examine this work closely.

The  1992 paper by Fenn and Rourke \cite{FR} also is very influential. At about that time, several of us were dabbling with the ideas of quandles (See for example \cite{onknots}), but there were only a few results that were being achieved.  A particularly, useful concept from the point of view of the current author was the notion of rack cohomology that was introduced\cite{FR2}. In works with other collaborators\cite{CJKLS},  the definition of rack cohomology was modified to that of quandle cohomology, and  we showed how to use quandle cocycles to produce invariants of classical knots and knotted surfaces in state-sum form. 

Brian Sanderson spoke about trunks in the 1995 Warsaw conference, and I was in the audience. I am pretty sure that I mentioned it to Masahico in passing, but neither of us recalls specifically such correspondence. The electronic record might produce a ``bit-trail," but for my part, I really did not think about the idea subsequently. 

Shortly after Laurel Langford finished her dissertation, the three of us (Laurel, Masahico, and I) met in February to discuss the possibility of using Neuchl's cocycles to define knotted surface invariants. Upon return from a meeting of the AMS in Manhattan, Kansas, Masahico and I reformulated the idea in terms quandle cocycles. As we began writing down the material, I began to recall Brian Sanderson's talk. I had some notes from the conference and the cocycle notion is clearly stated therein.  

In 1998, Seiichi Kamada and I spoke at length during the Knots in Hellas conference. Much of the pending discussion had to do with his upcoming travel arrangements to Mobile. He and Naoko arrived in Mobile in October of that year. Meanwhile, Dan Jelsovsky began pursuing his degree under Masahico's supervision. The original manuscript that contained the quandle cocycle invariant was posted \cite{CJKLS} on the arXiv in March 1999. 

A great deal of work by a number of authors in different contexts is based upon the quandle cocycle invariants and their generalizations. In the final section herein, I will summarize many of the results of which I know that use the cocycle invariants. 

Next, I turn to the basic definitions and examples.

\section{Fundamental Definitions and Results}
\label{Joyofsex}

\subsection*{Classical Quandle Theory}

\begin{definition}{\rm 
A {\it quandle}, $X$, is a set with a binary operation 
$(a, b) \mapsto a \lt b$
such that

\begin{enumerate}
\item[(I)]
 For any $a \in X$,
$a\lt a =a$.

\item[(II)] For any $a,b \in X$, there is a unique $c \in X$ such that 
$a= c\lt b$. 

\item[(III)] 
For any $a,b,c \in X$, we have
$ (a\lt b) \lt c=(a\lt c) \lt (b \lt c). $
\end{enumerate}

In Axiom (II), the element  $c$
that is uniquely determined from   given $a, b \in X$ such that $a=c\lt b$,
is denoted by $c=a \lt^{-1} b$. The reader may check that $\lt^{-1}$ also defines a quandle structure. 
A function $\phi: X \rightarrow  Y$ between quandles
or racks  is a {\it homomorphism}
if $(a \lt b)\phi = (a)\phi \lt  (b)\phi$ 
for any $a, b \in X$. In many contexts within the discussion herein, it will be convenient to write the function's name to the right of the argument, as above. Observe that in this notation, axioms (II) and (III) indicate that the quandle acts upon itself as a collection of {\it quandle automorphisms} --- one-to-one onto maps that preserve the quandle operation. The {\it automorphism group} of a quandle $X$ is the group ${\mbox{\rm Aut}}(X)$ of all bijective quandle homomorphisms. The {\it inner automorphism group}, $\Inn{(X)},$ is the subgroup of the automorphism group generated by $\{ \phi_z : z \in X \ \ \& \ \ (x)\phi_z= x \lt z \}.$

A quandle $X$ is said to be {\it connected} if for any elements $x,y \in X$ there is an inner automorphism that act transitively on $X$; thus there is an inner-automorphism $\phi \in \Inn{(X)}$ such that $(x)\phi=y$.   A quandle is said to be {\it homogeneous} if $\Aut{(X)}$ act transitively on $X$. If $X$ is connected, it is homogeneous.
A quandle is {\it involutory} if $\lt = \lt^{-1}$. 

} \end{definition}

\begin{example}{\rm 
The following are typical examples of quandles.  
A group $X=G$ with
$n$-fold conjugation
as the quandle operation: $a \lt b=b^{n} a b^{-n}$ or $a \lt b=b^{-n} a b^n$.
The notation Conj$(G)$ indicates the quandle defined for a group $G$ by 
$a \lt b=bab^{-1}$. 
Any subset of $G$ that is closed under such conjugation 
is also a quandle. 

Any $\Lambda (={\Z }[t, t^{-1}])$-module $M$  
is a quandle with
$a \lt b=ta+(1-t)b$, $a,b \in M$, 
that is 
called an {\it  Alexander  quandle}.
For a positive integer
$n$, the ring 
${\Z }_n[t, t^{-1}]/(h(t))$
is a quandle
for
a Laurent polynomial $h(t)$.
It 
is finite if the coefficients of the
highest and lowest degree terms
of $h$   are units in $\Z_n$.

Let $n$ be a positive integer, and 
for elements  $i, j \in \{ 0, 1, \ldots , n-1 \}$, define
$i \lt  j \equiv 2j-i \pmod{n}$.
Then $\lt $ defines a quandle
structure  called the {\it dihedral quandle},
  $R_n$.
This set can be identified with  the
set of reflections of a regular $n$-gon
  with conjugation
as the quandle operation, and it is also isomorphic to the Alexander quandle 
${\Z }_n[t, t^{-1}]/(t+1)$. 
As a set of reflections of the 
regular 
$n$-gon, 
 $R_n$ can be considered as a subquandle of ${\mbox{\rm Conj}}(\Sigma_n)$
where $\Sigma_n$ denotes the symmetric group on $n$ letters. }\end{example}

\begin{example}{\rm 
This example is due to Joyce \cite{Joyce} and Matveev\cite{Matveev}.

Let $G$ be a group, $H$ a subgroup, $s:G\rightarrow G$ an automorphism such that for each $h\in H$ $s(h)=h$. Define a binary operation $\lt_s=\lt$ on $G$ by $a\lt b=s(ab^{-1})b.$
Then $\lt$ defines a quandle structure on $G$. Axioms $I$ and $II$ are easily verified. For Axiom III, we have 
$$(a\lt b) \lt c= s( s(ab^{-1})b c^{-1})c= s^2(a)s^2(b^{-1})s(b)s(c)^{-1}c$$
while 
\begin{eqnarray*}\lefteqn{
(a \lt c) \lt (b \lt c)} \\ 
&= & s( s(ac^{-1})c (s(bc^{-1})c)^{-1} )s(bc^{-1})c \\ & =&  s^2(a)s^2(c^{-1})s(c)  [s^2(b)s^2(c^{-1})s(c)]^{-1} s(b) s(c^{-1})c\\
& = &  s^2(a)s^2(c^{-1})s(c)s(c)^{-1}[s^2(c^{-1})^{-1}]s^2(b)^{-1}s(b) s(c^{-1})c \\ & = &
s^2(a)s^2(b^{-1})s(b) s(c)^{-1}c. \end{eqnarray*}
This passes to a well-defined quandle structure on $G/H$ that is given by
$Ha \lt Hb= Hs(ab^{-1})b$. 
In particular, if $z\in Z(H)\cap H$ where $Z(H) = \{ z \in G: zh=hz \ {\mbox{\rm for all }} \  h\in H \}$, then
$Ha \lt Hb=Hab^{-1}zb$ defines a quandle structure. Let us denote the resulting quandle by $(G,H,z)$.}\end{example}

Let $X$ be any quandle. Let $F(X)$ denote the free group generated by $X$. Then $F(X)$
acts on $X$ as follows. 
 If $x\in  X$ also denotes the generator in $F(X)$ and $\cdot$ denotes the action of $F(X)$ on $X$, then $y \cdot x = y \lt x$;  furthermore 
$y \cdot x^{-1}= y\lt^{-1} x$. In general, a word $w= w_1w_2 \in F(X)$ acts as $y \cdot (w_1 w_2)= (y\cdot w_1)\cdot w_2$. This, in turn, can be written in terms of $\lt, \lt^{-1}$ and the generators in $X$. Thus there is a map from  $F(X)$ to a subgroup of the group of symmetries of $X$. Let $K$ denote the kernel of this map. Then $F(X)/K$ is isomorphic to a subgroup of the automorphism group of the quandle. In fact, since the action of $F(X)$ is induced by the quandle product, $F(X)/K$ is isomorphic to  $\Inn{(X)}$. 

Define the {\it enveloping group} of a quandle to be the group, 
$G_X$ that is generated by the elements of $X$ subject to the set of relations $a*b=bab^{-1}$ for all $a,b\in X$. The action of $F(X)$ also passes to an action of $G_X$ with kernel $K'$, and so we have in general $G_X/K'$ is isomorphic to a subgroup of $\Inn{(X)}$. 

Note that for a quandle $X$ we have have a map $\partial:X \rightarrow \Aut{(X)}$ given by $y \mapsto \phi_y$ where, $(x)\phi_y = x \lt y.$ Note also, that $\partial((x)\phi)= \phi^{-1} \partial(x) \phi.$ Thus a quandle 
is also a crossed $G$-set where $G=\Aut{(X)}$, but it has additional structure.

\begin{theorem}\label{thm:JM}{\rm \cite{Joyce,Matveev}}
Let $X$ denote a homogeneous quandle and suppose that $z$ is a fixed element of $X$. Consider $H \subset \Aut{(X)}$ to be the stabilizer of $z$. That is $H=\{ h \in \Aut{(X)}: (z)h=z \}$. Then there is a quandle isomorphism $X \cong (\Aut{(X)}, H, z)$.\end{theorem}
\noindent
{\it Proof.}
Let $x \in X$, and chose $\phi_x$ so that $x= (z)\phi_x$. Given $\phi_x$ and $\phi'_x$ that map $z$ to $x$, we see that $(z)\phi' \phi^{-1}=z$. So $\phi' \phi^{-1}\in H$. On the other hand, if $h\in H$, then $x= (z)\phi_x = zh\phi_x.$ Therefore, any element of  $H\phi_x$ maps $z$ to $x$. We have a function, $E:\Aut{(X)} \rightarrow X$, $E(\phi)=(z)\phi$. To see that $E$ is a quandle homomorphism, compute: 
$$E(\phi \lt \psi)= (z) ((\phi \psi^{-1}) \lt z)\psi= (z)(\phi \psi^{-1}) \psi \lt (z)\psi= E(\phi) \lt E(\psi).$$
The quandle map $E$ passes to  a bijection form $(\Aut{(X)}, H, \phi_z)$ to $X$. This completes the proof.

\begin{example}{\rm 
Let $QS_4$ denote the quandle, $\Z_2[t,t^{-1}]/(t^2+t+1)$. Observe that this example is the $4$-element field, but its quandle structure is given by
$a*b=ta +(1+t)b.$ Letting $[0]=0$, $[1]=1$, $[2]=t$, and $[3]=t+1$, we have the following quandle table:

\begin{center}
\begin{tabular}{|c||c|c|c|c||} \hline \hline
$R*C$ & $[0]$ & $[1]$ & $[2]$ &$[3]$ \\ \hline \hline
$[0]$ & $[0]$  & $[3]$ & $[1]$ & $[2]$ \\ \hline
$[1]$ & $[2]$ & $[1]$ &$[3]$ &$[0]$ \\ \hline
$[2]$ & $[3]$ & $[0]$ & $[2]$ & $[1]$ \\ \hline
$[3]$ & $[1]$ & $[2]$ & $[0]$ & $[3]$ \\ \hline \hline
\end{tabular}
\end{center}

Now $\Aut{(X)} \subset \Sigma_4=\Sigma_4(0,1,2,3)$, is generated by $(123)$, 
$(032)$, $(013)$, and $(021)$. So $\Aut{(X)}$ is the alternating group, $A_4= A_4(0,1,2,3).$
Let $H={\mbox{\rm Stab}}(\{0\})= \{ 1, (132),(123) \}$. The action of $[0]$ as an automorphism, is $\partial [0]= (123).$ The stabilizer of $0$ in $A_4$ is $H=\{1,(132),(123)\}$.
And the coset space is $\{ H, H(031), H(02)(13), H(013) \}$.}\end{example}

\begin{definition}{\rm 
 Suppose that
 $k:M^n \rightarrow \R^{n+2}$ is an embedding which is either smooth or  PL locally-flat. Let $N(k)$ denote an open tubular neighborhood of $k(M)$. 
The {\it fundamental quandle} $\pi_Q(k)$, of a codimension $2$ embedding is defined 
to be the set of homotopy classes of maps $\alpha: ([0,1], \{0,1\}, 1)\rightarrow (\R^{n+2}\setminus N(k), \partial N(k) \cup \{y_0\}, \{y_0\})$ where $y_0$ is a fixed base point pretty close to the boundary of the tubular neighborhood; later on, the base point will be chosen to be on the neighborhood. The homotopies between such maps are required to have their bottom boundaries on $\partial N(k)$ and their top boundaries fixed at the base point ($H(s,0) \in \partial N(k)$ while $H(s,1) = y_0$ for all $s\in [0,1]$, and $H(i,t)= \alpha_i(t)$ for $i=0,1$). If $\alpha$ and $\beta$ are such paths, then there is a unique oriented meridian $\mu_\beta$ that passes through the initial point of $\beta$. The quandle product is defined to be the path composition $\alpha * \beta = \alpha \beta^{-1} \mu_\beta \beta$. Here $\mu_\beta$ is the meridian that intersects the path $\beta$ at $\beta_0.$

A summary of the definition of the fundamental quandle is depicted in Fig.~\ref{funquan}.

\begin{figure}[htb]
\begin{center}
\includegraphics[width=4.5in]{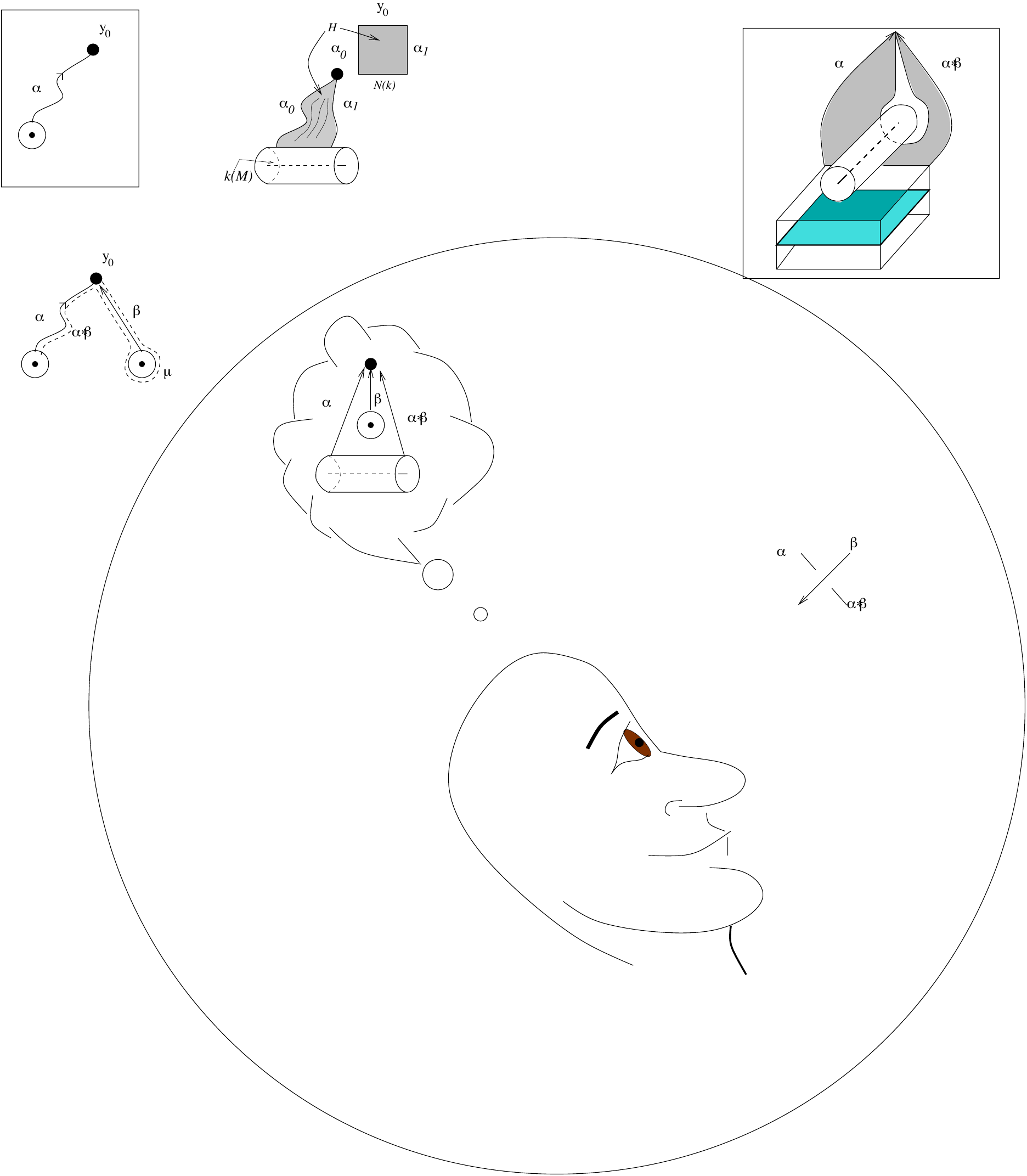}
\end{center}
\caption{The fundamental quandle} 
\label{funquan}
\end{figure}

}\end{definition}

Assume that $n=1$, so that the following discussion occurs in the classical case. 
The fundamental group $\pi_1(k)$ acts upon the fundamental quandle $\pi_Q(k)$ via path multiplication: given  representatives $\alpha: [0,1] \rightarrow M\setminus {\mbox{\rm int}} N(k)$ and $\gamma:([0,1], \{0,1\})\rightarrow ({\mbox{\rm int}} N(k), *)$
of the corresponding elements, then first travel along the quandle  path $\alpha$ from the tubular neighborhood to the base point and once around the loop $\gamma$. The result is a new path $[\alpha \cdot \gamma] \in \pi_Q.$

Assume that $k$ is knotted. Let  $G= \pi_1(S^3 \setminus k(S^1))$, and choose a base point on the boundary of the tubular neighborhood of $k(S^1)$. Then $G \subset \Aut{(\pi_Q(k))}$
except in the case of the unknot. (The unknot is singular since its fundamental quandle is trivial; hence the automorphism group is trivial while the fundamental group is $\Z$.) 
 Let $z \in \pi_Q$ denote the constant path, and let $H$ denote the peripheral subgroup that is induced by the inclusion of the torus into the knot complement. Then the elements of $H$  fix $z$ since the endpoint of $z$ can follow as it likes around a path on the tubular neighborhood that defines an element 
of $H$. See Fig.~\ref{action}. Moreover these are precisely the elements that fix $z$ because any element that fixes $z$ would result in the initial point of $z$ moving around the tubular neighborhood. So a loop that fixes $z$ is homotopic to a loop on the tubular neighborhood. Using the above quandle isomorphism, we have 

\begin{theorem} {\rm \cite{Joyce,Matveev}} The fundamental quandle is given as $\pi_Q(k)\cong(\pi_1(S^3 \setminus k(S^1)),  P, \mu_z)$ where $\mu_z$ is the meridonal element in the fundamental group that links $k(S^1)$ at the base point. \end{theorem}


\begin{figure}
\begin{center}
\includegraphics[width=3.5in]{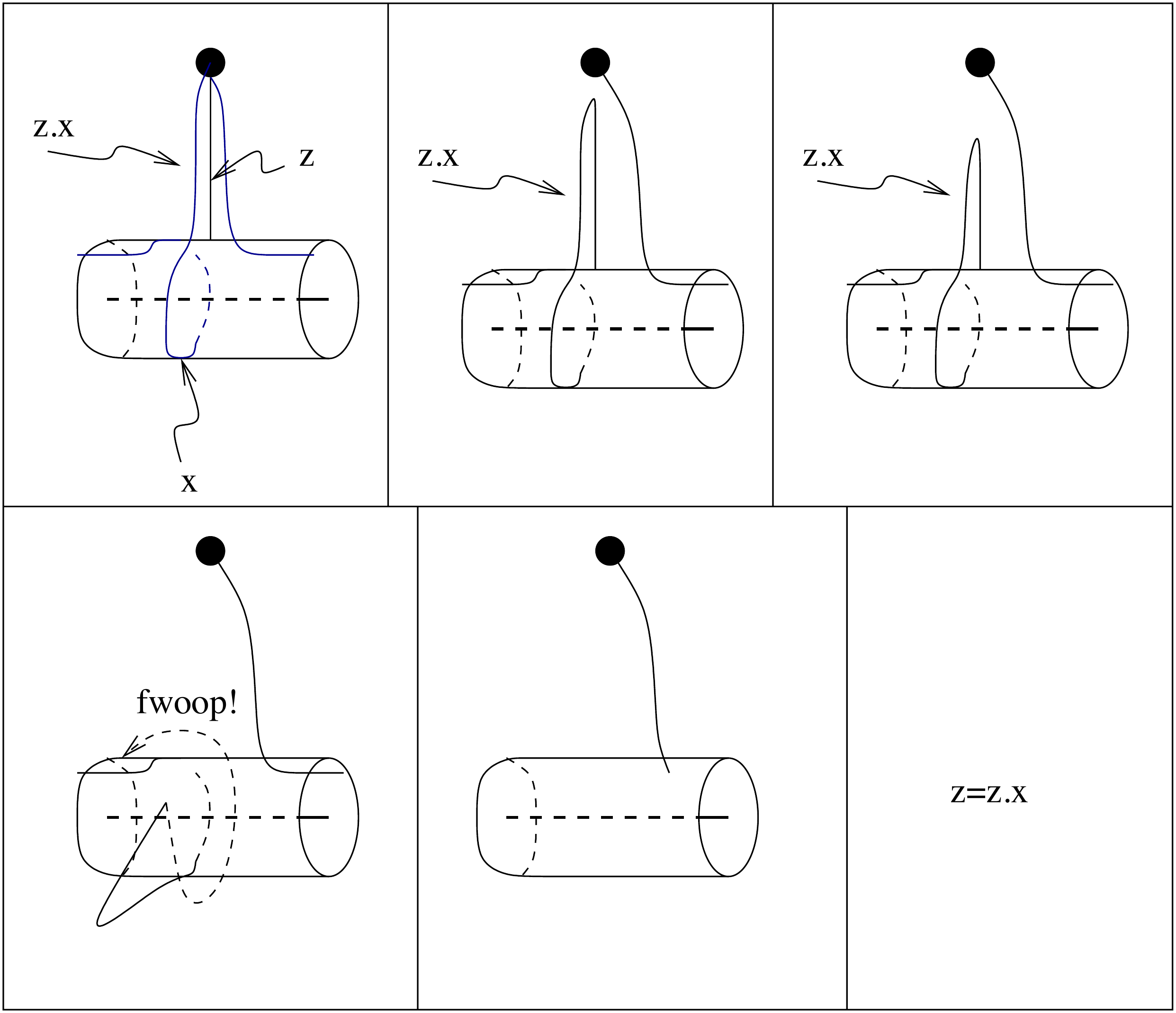}
\end{center}
\caption{The action of the peripheral subgroup}
\label{action}
\end{figure}

\subsection*{Quandle Colorings}

\begin{figure}
\begin{center}
\includegraphics[width=3.in]{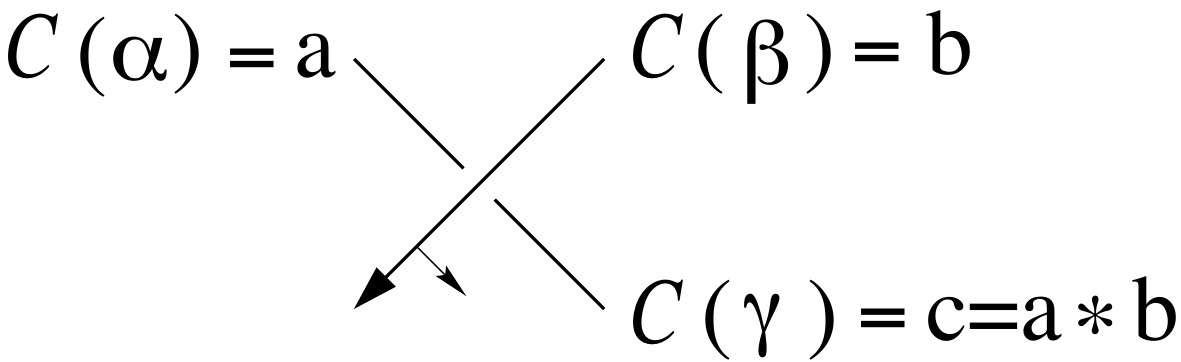}
\end{center}
\caption{Quandle coloring at a classical crossing}
\label{qcolor}
\end{figure}

\begin{definition}{\rm 
Let $X$ be a fixed quandle.
Let $K$ be a given oriented classical knot or link diagram,
and let ${\cal R}$ be the set of (over-)arcs. 
A vector perpendicular to an arc in the diagram 
is called  a {\it normal} or {\it normal vector}. We choose the normal so that 
the ordered pair  
(tangent, normal)  
 agrees with the orientation of the plane. This normal
  is called the {\it orientation normal}.  
A (quandle) {\it coloring} ${\cal C}$ is a map 
${\cal C} : {\cal R} \rightarrow X$ such that at every crossing,
the relation depicted in Fig.~\ref{qcolor} holds. 
More specifically, let $\beta$ be the over-arc at a crossing,
and $\alpha$, $\gamma$ be under-arcs such that the normal of the over-arc
points from $\alpha$ to $\gamma$.
(In this case, $\alpha$ is called the {\it source arc} and  $\gamma$ 
is called the {\it target arc}.) 
Then it is required that ${\cal C}(\gamma)={\cal C}(\alpha)\lt {\cal C}(\beta)$.
The color ${\cal C}(\gamma)$ depends only on the choice 
of orientation of the over-arc; therefore this rule defines the coloring
 at both positive and negative crossings.
The  colors ${\cal C}(\alpha)$, ${\cal C}(\beta)$
are called {\it source} colors.
}\end{definition}

The definition has a strict analogue for diagrams of knotted surfaces (and indeed higher dimensional knottings), but in order to present this, I will need to spend a little time on the notion of a knotted surface diagram. Specifically, suppose that $k:F \rightarrow \R^{4}$ is a smooth embedding of a closed surface (compact without boundary)
into $4$-dimensional space. The image may be perturbed slightly so that the image $f=p(k(F))$ under a projection $p:\R^4\rightarrow \R^3$ is a smooth surface in general position in $3$-space. The general position assumption means that each point
 $f(x)$  (for $x\in F$) has a neighborhood in $\R^3$ in which the image equivalent to any one of the following four scenarios: 
 \begin{enumerate}
 \item 
 $\{ (x,0,z) \in D^3:  x^2+ y=z^2 \le 1 \}$; here $D^3$ represents a $3$-dimensional ball of radius $1$; the point $a\in F$ that maps to $(0,0,0)$ is called a  {\it non-singular point};
 \item 
 $\{ (x,y,0)\in D^3:  x^2+ y^2 \le 1 \} \cup \{ (x,0,z) \in D^3:  x^2+ z^2 \le 1 \}$; the point $a \in F$ that maps to $(0,0,0)$ in this neighborhood is called a {\it double point};
 \item 
 $\{ (x,y,0)\in D^3:  x^2+ y^2 \le 1 \} \cup \{ (x,0,z) \in D^3:  x^2+ z^2 \le 1 \}\cup (0,y,z) \in D^3: y^2 +z^2 \le 1 \}$; the point $a\in F$ that maps to $(0,0,0)$ in this neighborhood is called a {\it triple point};
 \item The image of $f$ in a neighborhood of $p$ resembles the cone on a figure-$8$; the point $p$ in this case is called a {\it  branch point}. 
 \end{enumerate}

\begin{figure}
\begin{center}
\includegraphics[width=5in]{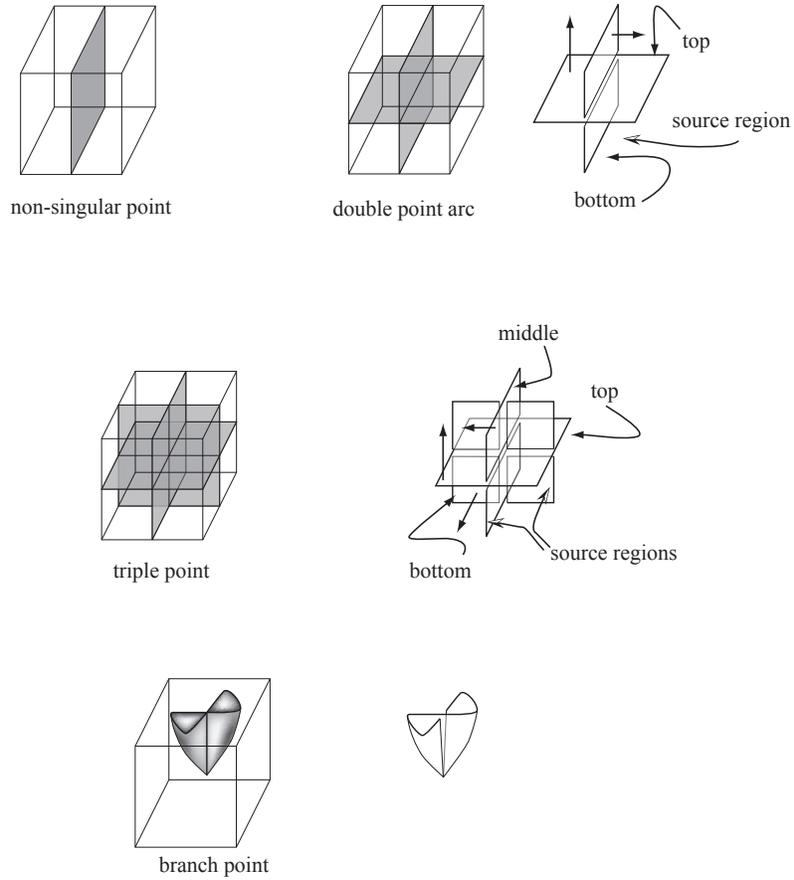}
\end{center}
\vspace{-1in}
\caption{Generic and broken surfaces}
\label{surfaces}
\end{figure}

The notion of a {\it broken surface diagram} is analogous to a knot diagram and can be constructed from such a generic projection. Specifically, at the singular points some of the sheets of the surface are more distant from the $3$-space into which the surface is being projected than the others are. Those that are farther away have a neighborhood of the double points removed within the diagrams. At a branch point, the differential of the function $f$ is singular; at a branch point a broken sheet converges to an unbroken sheet. Figure~\ref{surfaces} illustrates.

In the case that the surface $F$ is oriented, an orientation normal in $3$-space can be  chosen consistently across the broken sheets in the diagram. A {\it quandle coloring of a broken surface diagram} is defined as an assignment of quandle elements to each
sheet in a broken surface diagram such that the crossing rule depicted in Fig~\ref{qcolorsurf} is satisfied.  Specifically, suppose that two lower sheets are broken colored $a$ and $c$ and are separated by a sheet labeled $b$ with the normal to $b$ pointing from $a$ to $c$. Then $c=a\lt b$.

\begin{figure}
\begin{center}
\includegraphics[width=4.5in]{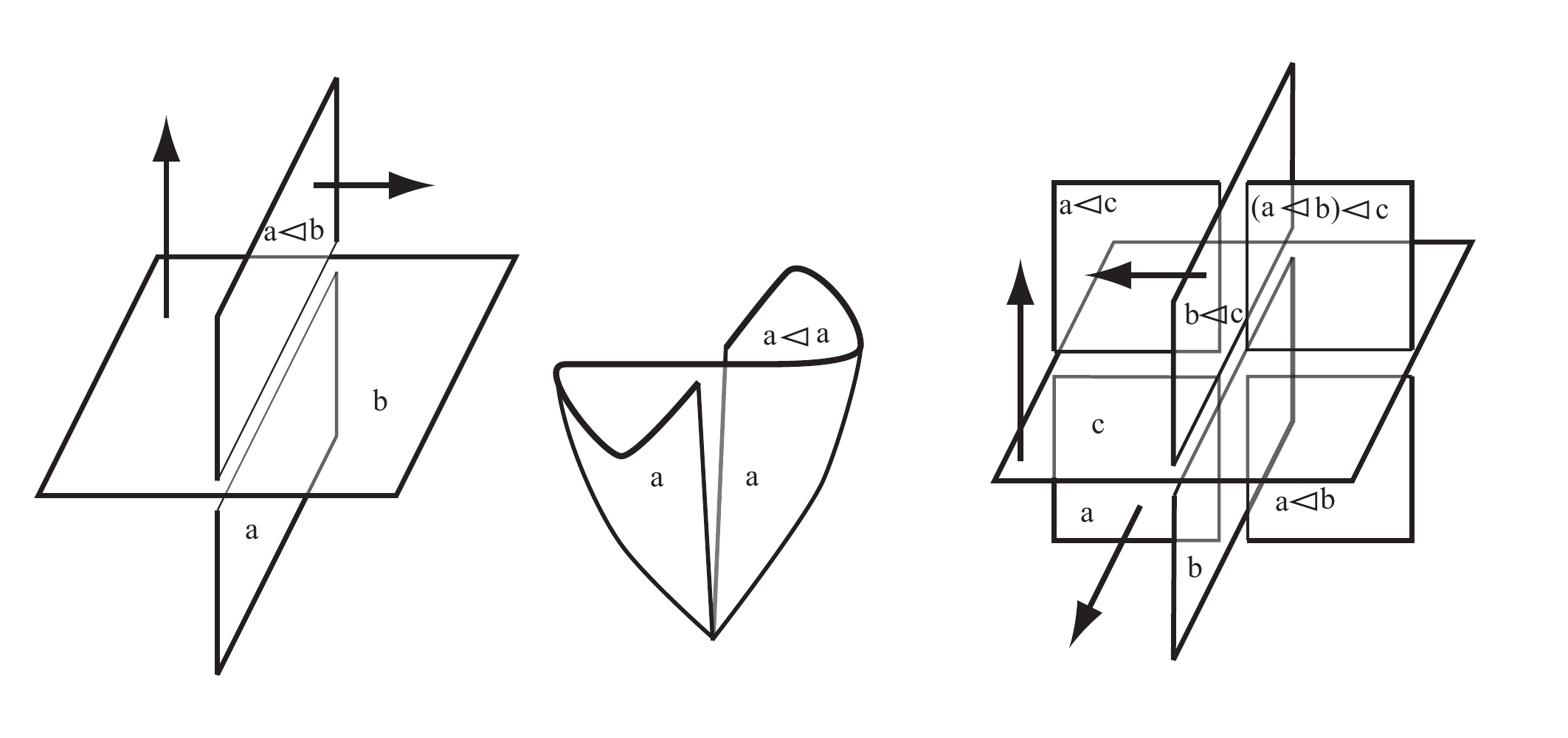}
\end{center}
\caption{Quandle  coloring of broken surfaces}
\label{action}
\end{figure}

In both the classical case and the knotted surface case, the quandle coloring is defined on a representative diagram of the given knot or knotted surface. If two classical diagrams differ by a Reidemeister move, and one of the diagrams has received a quandle coloring, then the quandle axioms allow us to color the second diagram in a unique fashion. In the knotted surface case, equivalent diagrams will differ by a finite sequence of Roseman moves \cite{Rose}. If two broken surface diagrams differ from each other by a Roseman move and one has a quandle coloring, then there is a uniquely induced quandle coloring on the other. In either case, we can envision the quandle coloring as a quandle homomorphism from the fundamental quandle of the codimension $2$ embedding to a test quandle. 

\subsection*{Quandle (Co)homology}

My collaborators and I presented\cite{CEGS} the discussion that follows to
 describe quandle cocycles and quandle homology in  the general setting as it was originally  introduced\cite{AG}. 

Let $X$ be a quandle. Let $\Omega (X)$ be the free ${\Z}$-algebra 
generated by 
$\eta_{x,y},$ $\tau_{x,y}$   for $x,y \in X$ such that 
$\eta_{x,y}$ is invertible for every $x,y \in X$.  
Define $\Z (X)$ to be the  quotient $\Z(X)= \Omega(X)/ R$ 
where $R$ is the ideal 
generated by 

\begin{enumerate}
\item \hfil $ \eta_{x\lt y,z}\eta_{x,y} -  \eta_{x \lt z,y \lt z}\eta_{x,z}$ \hfill \hfill
\item \hfil $ \eta_{x \lt y,z}\tau_{x,y} -  \tau_{x \lt z,y \lt z}\eta_{y,z}$ \hfill \hfill
\item \hfil $\tau_{x \lt y,z}- \eta_{x \lt z,y \lt z} \tau_{x,z}- 
 \tau_{x \lt z,y \lt z}\tau_{y,z}$ \hfill \hfill
\item \hfil $\tau_{x,x} + \eta_{x,x} -1$ \hfill \hfill
\end{enumerate}

The algebra  $\Z(X)$ thus defined is called the  {\it quandle algebra} 
over $X$. 
In $\Z(X)$, we define elements 
$\overline{\eta_{z,y}} =\eta^{-1}_{z\lt^{-1} y,y}$
and 
$\overline{\tau_{z,y}} = -\overline{\eta_{z,y}}\tau_{z\lt^{-1} y,y}$.
The convenience of such quantities will become apparent by examining type II 
moves.

A {\it representation} of $\Z(X)$ is an abelian group $G$ together with 
(1) a collection of automorphisms
$\eta_{x,y} \in {\mbox{\rm Aut}} (G)$, and (2) a collection of endomorphisms 
$\tau_{x,y}\in {\mbox{\rm End}} (G)$ such that the relations above hold. 
More precisely, there is an algebra homomorphism  
$\Z(X) \rightarrow {\mbox{\rm End}} (G)$, and we denote the 
image of the generators by the same symbols. Given a representation of 
$\Z(X)$ we say that $G$ is a $\Z(X)$-module, or {\it a quandle module}. 
The action of $\Z(X)$ on $G$ is written by the left action,  
and denoted by  $(\rho, g) \mapsto \rho g  (= \rho \cdot g = \rho (g) )$,
for  $g\in G$ and  $\rho \in End(G).$

\begin{example}{\rm \cite{AG}} \label{AGexample} {\rm 
Let $\Lambda = \Z[t,t^{-1}]$ denote the ring of Laurent polynomials. 
Then any  $\Lambda$-module $M$ is a $\Z(X)$-module for any quandle $X$,
by $\eta_{x,y} (a)=ta$ and $\tau_{x,y} (b) = (1-t) (b) $
for any $x,y \in X$.

Consider the  enveloping group  $G_X=\langle x \in X \ | \ x\lt y=yxy^{-1} \rangle$
(a.k.a. the {\it associated group} \cite{FR}). 
For any quandle $X$,
any $G_X$-module $M$  is a  $\Z(X)$-module by 
 $\eta_{x,y} (a)=y a$ and $\tau_{x,y} (b) = (1- x \lt y) (b) $,
where $x, y \in X$, $a, b \in M$.

} \end{example}

\begin{figure}[h]
\begin{center}
\includegraphics[width=3in]{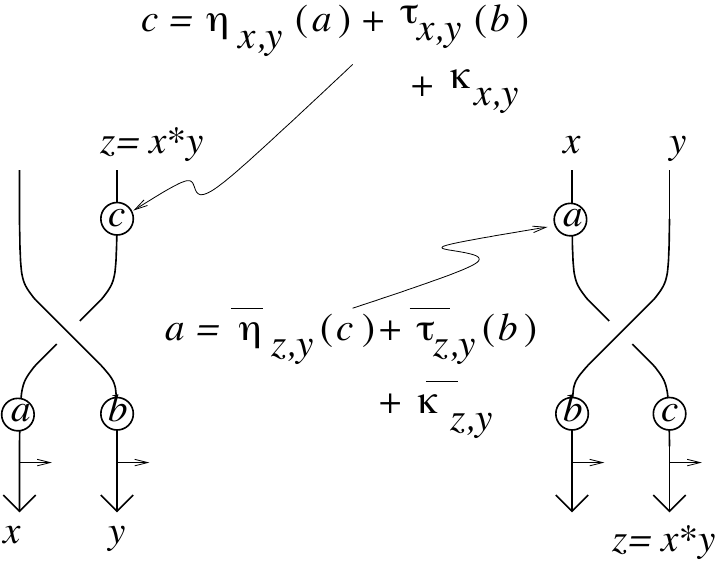}
\end{center}
\caption{The geometric notation at a crossing}
\label{gencolor} 
\end{figure}

\begin{figure}[h]
\begin{center}
\includegraphics[width=4in]{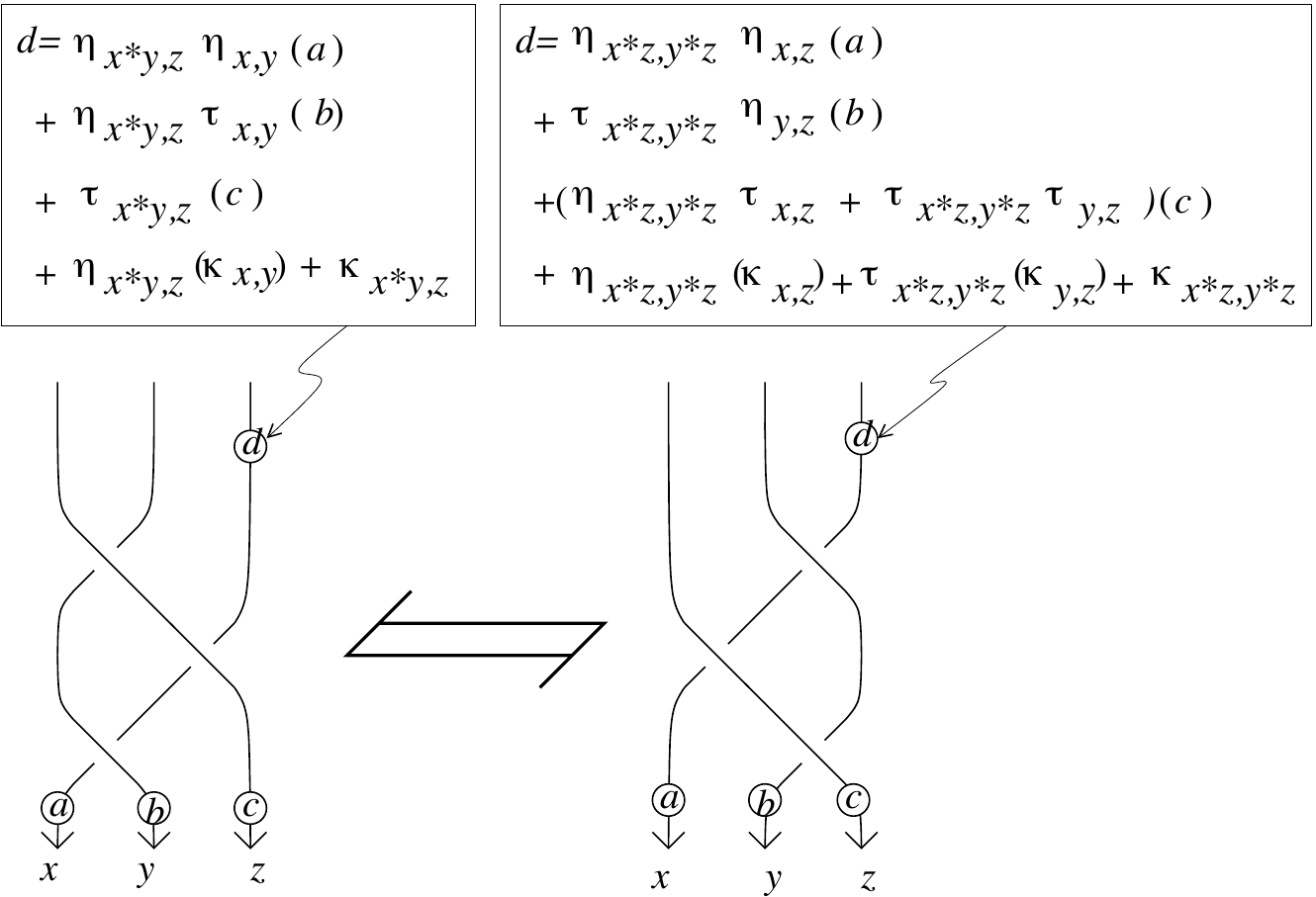}
\end{center}
\caption{Reidemeister moves and the quandle algebra definition}
\label{move3} 
\end{figure}

Figures~\ref{gencolor} and \ref{move3} 
indicate the geometric motivation for 
the quandle module axioms. 
For the time being, 
ignore the terms $\kappa_{x,y}$ in the figures.

Consider 
the free 
left 
$\Z ( X )$-module $C_n(X)=  \Z ( X ) X^n$ 
with basis $X^n$ (for $n=0$, $X^0$ is a singleton $\{ x_0 \}$, 
for a fixed element $x_0 \in X$). 
Boundary operators
 $\partial =\partial_n: C_{n+1}(X) \rightarrow C_n(X)$ 
are defined\cite{AG} by 
\begin{eqnarray*}\lefteqn{ \partial (x_1, \ldots, x_{n+1} )}\\ &=&
{\displaystyle (-1)^{n+1} \sum_{i=2}^{n+1} (-1)^i 
\eta_{ [x_1, \ldots, \widehat{x_i}, \ldots, x_{n+1} ], [x_i, \ldots, x_{n+1}] }
(x_1, \ldots, \widehat{x_i}, \ldots, x_{n+1} ) } \\
 & & \\
 & & - {\displaystyle (-1)^{n+1} \sum_{i=2}^{n+1} (-1)^i 
(x_1\lt x_i, \ldots, x_{i-1}\lt x_i, x_{i+1},  \ldots, x_{n+1} ) }  \\
 & & \\
 & & + (-1)^{n+1} 
 \tau_{[x_1, x_3, \ldots, x_{n+1}], [x_2, x_3, \ldots, x_{n+1}]}
(x_2, \ldots, x_{n+1} ) , \\[3mm]
\mbox{where} & & 
 [x_1, x_2, \ldots, x_n] = ( ( \cdots ( x_1 \lt  x_2 ) \lt  
x_3  ) \lt  \cdots ) \lt  x_n
\end{eqnarray*}
for $n > 0$, and 
$\partial_0(x) = - \tau_{x \bar \lt  x_0, x_0} x_0$ for $n=0$. 
In particular, the  
condition to be a 
$2$-cocycle 
for a 
$2$-cochain 
$\kappa:C_2(X) \rightarrow \Z(X)$ (or  ($\kappa: C_2(X) \rightarrow G) \in C^2(X;G)={\mbox{\rm Hom}}_{\Z(X)} (C_2(X),G)$ when a representation $\Z(X) \rightarrow G$ is fixed to define the twisted coefficients) in this homology theory is written as
$\kappa_{x,y}$ 
$$ \eta_{x\lt y, z} (\kappa_{x,y}) + \kappa_{x\lt y, z}
= \eta_{x\lt z, y\lt z} (\kappa_{x,z})  +  \tau_{x\lt z, y\lt z} (\kappa_{y,z})
 + \kappa_{x\lt z,y\lt z}, $$ 
in $\Z(X)$ (or in $G$) 
for any $x,y,z \in X$,
where $\kappa_{x,y}$ means $\kappa(x,y)$ for $(x,y)\in X^2.$ 
This is called a {\it generalized $($rack$)$ $2$-cocycle condition}.
When $\kappa $ further satisfies $\kappa_{x,x}=0$ for any $x \in X$,
 it is called a {\it generalized quandle $2$-cocycle.}

In case, $\eta_{x,y}=1$ and $\tau_{x,y}=0$ for all $x,y\in X$, then we recover the original quandle homology theory \cite{CJKLS}.

The (rack) homology group is defined to be the quotient of the kernel of $\partial_n$ modulo the image of $\partial_{n+1}$. 

\subsection*{Dynamical cocycles}
Let $X$ be a quandle and $S$ be a non-empty set. 
Let $\alpha: X \times X \rightarrow 
\mbox{\rm Fun}(S \times S, S)=S^{S \times S}$ be a function,
so that for $x,y  \in X$ and $a, b \in S$ we have 
$\alpha_{x,y} (a,b) \in S$. 

Then it is checked by computations 
that $S \times X$ is a quandle by the operation
$(a, x)\lt (b, y)=(\alpha_{x, y}(a,b )  , x\lt y )$,
where $ x\lt y$ denotes the quandle 
product in $X$, 
if and only if $\alpha$ satisfies the following conditions:

\begin{enumerate}
\item $\quad\alpha_{x, x} (a,a)= a$ for all $x \in X$ and 
$a \in S$; 
\item $\quad\alpha_{x,y} (-,b): S \rightarrow S$ is a bijection for 
all $x, y \in X$ and for all $b \in S$;
\item $\quad 
\alpha_{x \lt  y, z} ( \alpha_{x , y}(a, b), c)
=\alpha_{x \lt  z, y\lt z}(\alpha_{x,z}(a,c), \alpha_{y,z}(b,c) )$
for all $x, y, z \in X$ and $a,b,c \in S$. 

\end{enumerate}

Such a function $\alpha$ is called a {\it dynamical quandle cocycle}
\cite{AG}.
The quandle constructed above is denoted by $S \times_{\alpha} X$, 
and is called the {\it extension} of $X$ by a dynamical cocycle $\alpha$.
The construction is general, as Andruskiewitsch and Gra\~{n}a
 show:

\begin{lemma}{\rm \cite{AG}} \label{AGlemma} 
Let $p: Y \rightarrow X$ be a surjective quandle homomorphism 
between finite quandles 
such that the cardinality of $p^{-1}(x)$ is a constant for all $x \in X$. 
Then $Y$ is isomorphic to an extension  $S \times_{\alpha} X$ of $X$ 
by some dynamical cocycle on 
the
set $S$ 
such that $|S|=|p^{-1}(x)|$.
\end{lemma}

\begin{example}{\rm {\cite{AG}}} {\label{whathesaid}} 
{\rm
Let $G$ be  a $\Z(X)$-module for the quandle $X$,
and $\kappa$ be a generalized $2$-cocycle. 
For $a,b \in G$, let 
$$\alpha_{x,y}(a,b) =  \eta_{x,y}(a)  +  \tau_{x,y}(b)
+ \kappa_{x,y} .$$ 
Then it can be verified directly that $\alpha$ is a dynamical cocycle. 
In particular, even with $\kappa=0$, 
a $\Z(X)$-module structure on the abelian 
group $G$ defines a quandle structure
$G\times_\alpha X$. }\end{example}


We  \cite{CEGS}  used these definitions of quandle cocycles to define invariants of knots and knotted surfaces in the case that the quandle module $M$ was a module over the  enveloping group $G_X$ of the quandle, and the actions of $\eta$ and $\tau$ were given as:
 $\eta_{x,y} (a)=y a$ and $\tau_{x,y} (b) = (1- x\lt y) (b) $,
where $x, y \in X$, $a, b \in M$. As an important application, we showed that a large number of knotted $2$-spheres are non-invertible. Not only that, but they remain non-invertible even when trivial $1$-handles are attached. 

I think that the best plan to explain the quandle cocycle invariants is to give the original definition (in the case $\eta_{x,y}=1$ and $\tau_{x,y}=0$ for all $x,y \in X$), and leave the curious to read the original sources.

\subsection*{Cocycle invariants of classical knots and knotted surfaces}

Let a quandle $X$ be given, and assume that $\phi: X \times X \rightarrow A$ is a quandle $2$-cocycle that takes values in an Abelian group $A$. 
Thus $$\phi(a,b)+ \phi(a\lt b , c) = \phi(a,c) + \phi( a\lt c , b \lt c),$$and 
$$\phi(a,a)=0$$
for all $a,b,c \in X.$

Consider a knot diagram ${\mathcal  D}(k)$ of a knot $k:S^1 \rightarrow \R^3$. Choose a quandle coloring of  ${\mathcal  D}(k)$ by $X$. In a neighborhood of a crossing, consider the under-arc away from which the normal to the over-arc points. The quandle color, $a$, on this arc is called the {\it source color of the under-arc}.  The color, $b$, on the over-arc is its source color. Consider the value of the cocycle $\pm \phi(a,b)$; this is called the {\it Boltzmann weight} of the crossing. The sign coincides with the sing of the crossing. The weight of the coloring is the sum (in $A$) over all crossings of the diagram ${\mathcal  D}(k)$: {\it i.e.} $B(\phi, C) = \sum_{\rm Crossings} \pm \phi(a,b)$
Then the {\it quandle cocycle invariant} of $k$ is the multiset of values  $\{ B(\phi, C) \}_C$ taken over all colorings $C$ of the diagram.

In the case of knotted surfaces, the definition of the cocycle invariant is similar.
The diagram of a triple point is considered. There are  four sheets that are labeled as under sheets,  two sheets that are labeled middle sheets, and one sheet that is labeled a top sheet. The normal to the top sheet points away from one of the middle sheets; that middle sheet is a {\it source  region}. The normal to both the top and middle point away from one of the four bottom sheet; this is also a {\it source region}.

We assume that the knotted surface diagram has been colored by a test quandle $X$. Let the color of the source region to the bottom be $a$, let  that of the source region in the middle be $b$, and let that of the color of the top sheet be $c$. The {\it sign of the triple point} is the orientation class of the triple (normal to top, normal to middle, normal to bottom) in comparision to the right handed orientation of $3$-space.

A quandle $3$-cocycle $\theta: X \times X \times X \rightarrow A$ 
is a function that takes values in an Abelian group $A$ and that satisfies:
$$\theta(a,b,c) + \theta(a\lt c, b\lt c, d) + \theta(a,c,d) = \theta(a,b,d)+ \theta(a\lt b, c, d) + \theta(a\lt d, b\lt d ,  c\lt d)$$and 
$$\theta(a,a,b)=\theta(a,b,b) = 0$$

At each triple point, the cocycle is evaluated on the ordered triple:  (the color on the source region of the lower sheet, the color on the source region of the middle sheet, and the color on top sheet). The values of the cocycles, $\pm \theta(a,b,c)$ are summed over the set of triple points where the sign is chosen to coincide with the sign of the triple point. And as in the classical case, the multiset of values over all possible colorings is the value of the cocycle invariant on the knotted surface. 

In both the classical case and the knotted surface case, the cocycle conditions, and the sign choices are exactly what are needed in order to show that the value of the invariant does not depend upon the choice of diagram representing the knotting. 

\bigskip

\noindent {\bf Some non-trivial cocycles.} Unfortunately, we do not have a lot of high-tech methods (such as projective resolutions or applications of homological algebra) to produce non-trivial cocycles. Some early non-trivial cocycles were determined by means of computer calculations. One of these was generalized by Mochizuki \cite{Mochi} via methods from algebraic geometry. Newer results will be summarized in the last section of the paper, but it is worth mentioning here that in the case of dihedral quandles, Niebrzydowski and Przytycki \cite{NP,NP2,NP3} have developed methods and conjectures about the higher order homology in the dihedral cases.

\begin{example} {\rm For  the Alexander quandle $QS_4= \Z[t]/(2,t^2+t+1)$, the  function, $$\phi(x,y)=  \left\{ \begin{array}{lr} 0 & {\mbox{\rm if}} \ \  x=[3], y= [3], \ \ {\mbox{\rm or}} \ \ x=y \\
1 & {\mbox{\rm otherwise}} \end{array} \right.$$
that takes values in $\Z_2$ is a quandle $2$-cocycle. }\end{example}

\begin{example}{\rm In the case of the dihedral quandle $R_p=\{0,1,2, \ldots, p-1 \}$ of order $p$ when $p$ is a prime ($a \lt b = 2b -a \pmod{p}$), a generator of $H^3_{\rm Q}(R_p, \Z_p) \cong \Z_p$ is represented by the function:
$$\theta_p(x,y,z) =(x-y)\frac{(2z-y)^p + y^p - 2 z^p}{p}.$$
The cocycle $\theta_p$ is Shin Satoh's simplification of Mochizuki's original construction. }\end{example}

\subsection*{Quandles with Good Involutions}

An interesting class of quandles consist of quandle that have a good involution. Let me start with an example. Consider a subset $Q$ of a group that is closed under conjugation and inversion. For elements $a,b \in Q$ define $a\lt b= ba b^{-1}$, and define $\rho(a)=a^{-1}$. The map $\rho$ is an involution, but not necessarily a quandle map. We have that $\rho(a \lt b)= \rho(a) \lt b$, and $a \lt \rho(b) = a \lt^{-1} b$. In general, an involution  $\rho$ on a quandle $X$ that satisfies these two properties will be called {\it a good involution,} and if $X$ possesses such an involution, then $(X, \rho)$ will be said to be {\it a symmetric quandle.}  The material of the current section is a summarization of material that originally appeared in our paper\cite{COS}.

The associated group, $G_{(X,\rho)}$ of a symmetric quandle $(X, \rho)$ was defined\cite{K, KO}
 by
$G_{(X, \rho)} = \langle x \in X : x \lt y= y^{-1} x y, \ 
 \rho(x)=x^{-1} \rangle$.    
For a symmetric quandle $(X, \rho)$, an $(X, \rho)$-set is a set $Y$ 
equipped with a right action of the
associated group $G_{(X, \rho)} $.
The notation $yg$ or $y \cdot g$ indicates the image of an element $y \in Y$  
under 
the action 
of  
$g \in G_{(X, \rho)} $. 
The   
 following 
three 
formulas hold:
$y \cdot (x_1 x_2) = (y \cdot x_1) \cdot x_2, \ $  
$y \cdot (x_1 \lt x_2) = y \cdot (x_2^{-1} x_1 x_2), \ $ and
$y \cdot (\rho( x_1 ) ) = y \cdot (x_1^{-1} ) $, 
for $x_1, x_2 \in X$ and $y\in Y$. 

In the case of a symmetric quandle, the homology theory of a quandle may be modified as follows to construct a symmetric homology theory.

\begin{sloppypar}
Let $Y$ be an $(X, \rho)$-set 
which may be empty. 
Let $C_{n}(X)_Y$ be the free abelian group generated by 
$(y, x_1,\ldots,x_n)$, where $y \in Y$ and $x_1, \ldots, x_n \in X$. 
For a positive integer $n$, 
let 
$C_0(X)_Y=\Z (Y)$, the free abelian group 
generated by  
$Y$, and 
set $C_n(X)_Y=0$ otherwise. 
(If $Y$ is empty, then  define $C_0(X) = 0$).  
Define the boundary 
homomorphism
\mbox{$\partial_n: C_{n}(X)_Y \longrightarrow C_{n-1}(X)_Y$}
by
\begin{eqnarray*}
\lefteqn{\partial_{n}(y, x_1,\ldots,x_n)
=\sum_{i=1}^{n}(-1)^i[(y, x_1,x_2,\dots,x_{i-1},\hat{x_i}, x_{i+1},\ldots,x_{n})}\\
&-&(y \cdot x_i, x_1\lt x_i,x_2\lt x_i,\ldots,x_{i-1}\lt x_i,\hat{x_i},  x_{i+1},\ldots, x_{n})]
\end{eqnarray*}
for $n\geq 1$ and $\partial_{n}=0$ for $n\leq 1$. Then
$C_{*}(X)_Y =\{C_{n}(X)_Y,\partial_n\}$ is a chain complex \cite{FRS}. 
Let
$D_{n}^Q(X)_Y$ be the subgroup 
of $C_{n}(X)_Y$ generated by 
$\cup_{i=1}^{n-1} \{ (y, x_1,\ldots, x_n )  \ | \ x_i=x_{i+1} \}$, and 
let 
$D_{n}^\rho (X)_Y$ be the subgroup 
of $C_{n}(X)_Y$ 
generated by $n$-tuples of  the form
$$ (y, x_1,\ldots, x_n ) + (y \cdot x_i , x_1\lt x_i ,\ldots, x_{i-1}\lt x_i , \rho(x_i), x_{i+1}, \ldots, x_n )$$     where $ y \in Y, \ \ $  $x_1, \ldots, x_n \in X, \ \ $ and $ i\in \{ 1, \ldots, n-1\}.$ 
Then $\{ D_{n}^Q(X)_Y, \partial_n \}$ and $\{ D_{n}^\rho (X)_Y, \partial_n \}$ are 
subcomplexes of $C_n$ \cite{KO}, and chain complexes 
$C_*^{R} (X)_Y$, $C_*^{Q} (X)_Y$, $C_*^{R, \rho} (X)_Y$, $C_*^{Q, \rho} (X)_Y$
are defined, respectively, from chain groups 
$C_n^{R} (X)_Y=C_n(X)_Y$, $C_n^{Q} (X)_Y=C_n (X)_Y/ D_n^Q (X)_Y$, 
$C_n^{R, \rho} (X)_Y=C_n (X)_Y/ D_n^\rho (X)_Y$, 
$C_n^{Q, \rho} (X)_Y = C_n (X)_Y / (D_n^Q  (X)_Y+ D_n^\rho (X)_Y)$. 
Their respective homology groups \cite{KO} are denoted by
$H_*^{R} (X)_Y$,  $H_*^{Q} (X)_Y$, $H_*^{R, \rho} (X)_Y$, 
and  $H_*^{Q, \rho} (X)_Y$, respectively. 
When $Y=\emptyset$, this subscript is dropped. 
Corresponding cohomology groups are defined as usual, as well as 
(co)homology groups with other coefficient groups.
\end{sloppypar}

Let a 
 surface diagram $D$ of a surface-knot $F$ be given. We cut the diagram further into {\it semi-sheets} by considering the upper sheets 
 also 
 to  
 be broken along the double point arcs. Observe that in the local picture of a branch point there are two semi-sheets, at a double point there are $4$ semi-sheets, and at a triple point, there are $12$ semi-sheets.

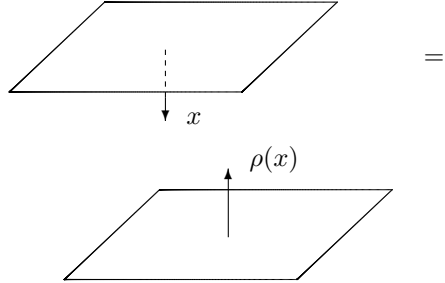
\begin{figure}[htb]
\begin{center}
\begin{minipage}{132pt}
\begin{picture}(132,70)

\qbezier(7,22)(7,22)(95,22)
\qbezier(7,22)(7,22)(43,56)
\qbezier(131,56)(43,56)(43,56)
\qbezier(131,56)(95,22)(95,22)

\put(66,20){\vector(0,-1){9}}
\multiput(66,38)(0,-4){5}{\line(0,-2){2}}
\put(74,10){$x$}

\end{picture}
\end{minipage}
\hspace{1cm}$=$\hspace{1cm}
\begin{minipage}{130pt}
\begin{picture}(130,70)

\qbezier(7,22)(7,22)(95,22)
\qbezier(7,22)(7,22)(43,56)
\qbezier(131,56)(43,56)(43,56)
\qbezier(131,56)(95,22)(95,22)

\put(69,38){\vector(0,1){26}}
\put(77,65){$\rho(x)$}

\end{picture}
\end{minipage}

\end{center}
\caption{A basic inversion}
\label{inversion}
\end{figure}

Let $(X,\rho)$ 
denote a 
symmetric 
quandle, and  let $Y$ 
denote 
an $(X, \rho)$-set.  
The surface diagram $D$ 
has 
elements of $X$ assigned 
to the sheets and elements of $Y$ assigned to regions separated by 
the projection such that the following conditions are satisfied.
\begin{itemize}
\item (Quandle coloring rule on over-sheets) 
Suppose that two adjacent semi-sheets coming from an over-sheet of $D$ about a double curve are labeled by $x_1$ and $x_2$. If the normal orientations are coherent, then $x_1=x_2$, otherwise $x_1= \rho(x_2)$. 
\item (Quandle coloring rule on under-sheets) 
Suppose that two adjacent 
under-sheets
$e_1$ and $e_2$ 
are separated along a 
double curve 
and 
are labeled by $x_1$ and $x_2$. 
Suppose that one of the two semi-sheets
coming from the over sheet 
of $D$, say $e_3$, is labeled by $x_3$.  We assume that 
a local 
normal orientation of $e_3$ 
points 
from $e_1$ to $e_2$.  If the normal orientations of $e_1$ and $e_2$ are coherent, then $x_1\lt x_3 = x_2$, otherwise $x_1\lt x_3= \rho(x_2)$. 
\item (Region colors)
Suppose that two adjacent regions $r_1$ and $r_2$ which are separated by a semi-sheet, say $e$, are
labeled by $y_1$ and $y_2$, where $y_1, y_2 \in Y$. Suppose that the semi-sheet $e$ is labeled by $x$.
 If the normal orientation
of $e$ points from $r_1$ to $r_2$,  then $y_1 \cdot x = y_2$. 
\item
 An equivalence relation (of 
a local normal orientation assigned to each semi-sheet and a quandle element associated to this local orientation) is generated by the following rule 
({\it basic inversion}\/):  
Reverse 
the normal orientation of a semi-sheet and 
change 
the element $x$ assigned the sheet by $\rho(x)$. 
\end{itemize}
A {\it symmetric quandle coloring}, 
 or an {\it $(X, \rho)_Y$-coloring, }
 of a surface-knot diagram is such an equivalence class of symmetric quandle colorings. 
 See Fig.~\ref{inversion}.

\begin{figure}[t]
\begin{center}
\begin{minipage}{130pt}
\begin{picture}(130,130)

\put(60,75){\setlength{\fboxrule}{18pt}\fcolorbox[gray]{.9}{.9}{}}
\thicklines
\thinlines

\qbezier(7,47)(7,47)(95,47)
\qbezier(7,47)(7,47)(29,68)
\qbezier(131,81)(95,47)(95,47)
\qbezier(131,81)(131,81)(110,81)

\qbezier(50,48)(50,48)(50,90)
\qbezier(84,122)(50,90)(50,90)
\qbezier(84,122)(84,122)(84,107)
\qbezier(50,48)(50,48)(69,66)

\qbezier(50,8)(50,8)(50,43)
\qbezier(50,8)(50,8)(69,26)
\qbezier(54,47)(50,43)(50,43)

\qbezier(71,68)(71,68)(110,68)
\qbezier(71,68)(71,68)(71,107)
\qbezier(71,107)(71,107)(110,107)
\qbezier(110,68)(110,68)(110,107)

\qbezier(71,28)(71,28)(110,28)
\qbezier(71,28)(71,28)(71,47)
\qbezier(110,28)(110,28)(110,61)

\qbezier(30,28)(30,28)(50,28)
\qbezier(30,28)(30,28)(30,47)

\qbezier(30,68)(30,68)(30,107)
\qbezier(30,107)(30,107)(66,107)
\qbezier(30,68)(30,68)(50,68)
\qbezier(66,106)(66,107)(66,107)

\put(34,86){\footnotesize\colorbox {white}{\textcolor{black}{$x_2$}} }

\put(50,80){\vector(-1,0){8}}
\multiput(55,80)(4,0){3}{\line(-1,0){2.0}}

\put(99,107){\vector(1,1){5}}

\qbezier(87,95)(87,95)(89,97)
\qbezier(91,99)(91,99)(93,101)
\qbezier(95,103)(95,103)(97,105)
\qbezier(99,107)(99,107)(101,109)

\put(80,47){\vector(0,-1){8}}
\multiput(80,50)(0,4){3}{\line(0,-1){2.0}}

\put(78,80){$[y]$}

\put(83,36){\colorbox [gray]{1.0}{\textcolor{black}{ \footnotesize$x_3$}}}
\put(100,118){ \footnotesize$x_1$}
\put(40,-10){ 
{ $+(y, x_1, x_2, x_3)$}}

\end{picture}
\end{minipage}

\end{center}
\caption{The weight of a triple point}
\label{weight}
\end{figure}

Let $(D,C_D)$ be an $(X,\rho)_Y$-colored diagram of an $(X,\rho)_Y$-colored surface-link $(F,C)$.
For a triple point of $D$, define the {\it weight} as follows:
Choose one of eight $3$-dimensional regions around the triple point and call the region a {\it specified region}.
There exist 12 semi-sheets around the triple points.
Let $S_T$, $S_M$ and $S_B$ be the three of them that face the specified region, where $S_T$, $S_M$ and $S_B$ are in the top sheet, the middle sheet and the bottom sheet at the triple point, respectively.
Let $n_T$, $n_M$ and $n_B$ be the normal orientations of $S_T$, $S_M$ and $S_B$ which 
point away from the specified region. 
Consider a representative of $C_D$ such that the normal orientations of $S_T$, $S_M$ and $S_B$ are $n_T$, $n_M$ and $n_B$ and
let $x_1$, $x_2$ and $x_3$ be the labels assigned the semi-sheets $S_B$, $S_M$ and $S_T$,  with the normal orientations  $n_B$, $n_M$ and $n_T$,  respectively. 
Let $y$ be the label assigned to the specified region. 
The weight is  
$\epsilon (y, x_1, x_2, x_3)$, where $\epsilon $ is $+1$ (or $-1$) if the triple of the normal orientations $(n_T,n_M,n_B)$ does (or does not, respectively) 
match 
the orientation of $\mathbb R^3$.
See Figure~\ref{weight}.
The sum $\sum_{\tau} \epsilon (y, x_1, x_2, x_3)$ of the weight over all 
triple points of a colored diagram $(D,C_D)$ represents a $3$-cycle $[c_D] \in C_3^{Q, \rho} (X)_Y$
\cite{KO}. 
Colored diagrams and cycles represented by colored diagrams are 
similarly defined 
when region colors are absent, or equivalently, when $Y=\emptyset$.

We \cite{COS} demonstrated (among other things) that for each positive integer $N$ there is a non-orientable surface embedded in a $3$-manifold $\times [0,1]$ such that the projection to the $3$-manifold is a generic surface with more than $N$ triple points. 

In addition, we produced the following family of quandles: 

\begin{theorem}{\rm {\cite{COS}} }\label{extthm}
For each positive integer $n$, there is an extension $\tilde{R}_{2n+1}$ 
of $R_{2n+1}$ with a non-trivial good involution $\rho$   
that is 
connected and 
is not involutory. 
\end{theorem}

Let $e_j$ denote the  column vector  in $\R^m$  whose $j\/$th entry is $1$ and the remaining entries are  each $0$. A {\it signed permutation matrix} is a 
square matrix of size $m$ 
matrix  whose columns are 
of the form $(\pm e_{\sigma (1)}, \pm e_{\sigma (2)},\ldots, \pm e_{\sigma (m)})$ where $\sigma \in \Sigma_m$ is a permutation. The set of signed permutation matrices form a group $H_m$ of order $2^m m!$ that is called the {\it hyper-octahedral group}.  
Define the group $SH_m$ to be 
the signed permutation matrices of determinant $1$.

 To avoid extra subscripts, we write $(\pm e_{\sigma (1)}, \pm e_{\sigma (2)},\ldots, \pm e_{\sigma (m)})$ as $(\pm {\sigma (1)}, \pm {\sigma (2)},\ldots, \pm {\sigma (m)}).$ 
 Then the matrix multiplication, in this notation, is written  as 
 \begin{eqnarray*} \lefteqn{(\epsilon_1 \cdot {\sigma (1)},\  \ldots,\  \epsilon_m \cdot {\sigma (m)})\cdot  (\delta_1 \cdot {\tau (1)}, \ \ldots, \ \delta_m \cdot {\tau (m)})} \\
& = & 
(\epsilon_{\tau(1)} \delta_1\cdot \sigma (\tau(1)),  \ \ldots, \ \epsilon_{\tau(m)} \delta_m\cdot \sigma (\tau(m)))
,\end{eqnarray*}
where $\epsilon_i=\pm 1$ and $\delta_j=\pm 1$ for $i, j= 1, \ldots, m$. 
For example, 
$(1,5,4,-3,-2)\cdot(5,1,2,3,4)=(-2,1,5,4,-3)$ 
while 
$(5,1,2,3,4)\cdot(1,5,4,-3,-2)=(5,4,3,-2,-1).$

We identify a subgroup of $SH_m$ that maps onto the dihedral group.  
Let $m=2n+1$. Consider the subgroup $G_{2n+1}$   of $SH_{2n+1}$ that is generated by the pair of elements $a=(1, 2n+1, 2n, \ldots, n+2, -(n+1), \ldots, -3, -2)$ and $b=(2n+1,1,2,\ldots , 2n)$.

The dihedral group, 
$D_{2(2n+1)}$, will be identified  
with the image of its faithful representation 
in permutation matrices.
Specifically, 
we 
identify 
the reflection 
$x$ with $(1, 2n+1, 2n, \ldots, 2)$, 
the rotation 
$y$ with $(2n+1,1,2,\ldots , 2n)$, 
and $D_{2(2n+1)}$ with the subgroup of permutation matrices generated by 
these two elements. 
Then the map that takes each matrix $(a_{ij})$ to 
$(|a_{ij}|)$ defines  
a group homomorphism onto the dihedral group: 
$f: G_{2n+1} \rightarrow D_{2(2n+1)}$, 
such that $f(a)=x$ and $f(b)=y$.

\begin{lemma}{\rm \cite{COS}} \label{orderlem} The order of $G_{2n+1}$ is $(2n+1)\cdot 2^{2n+1}.$ 
The centralizer of $a$, $C(a)=\{ c \in G_{2n+1}: ac=ca\}$,  is generated by the elements 
$$(1,\;  \epsilon_2 \;(2n+1),\;  \epsilon_3 \ (2n), \; \ldots, \; \epsilon_{n+1} \ (n+2), \; - \epsilon_{n+1} \ (n+1), \; \ldots ,\;  -\epsilon_2\ (2))$$ where $\epsilon_j = \pm 1$ for  $j=2,\ldots , n+1.$ In particular, 
$|C(a)|= 2^{n+1}.$ \end{lemma}
Our paper\cite{COS} contains the proof.

 Let $H=C(a)$ denote the centralizer of $a$. Consider the quandle structure $\tilde{R}_{2n+1}= (G_{2n+1}, H, a)$  given by $Hu \triangleleft  Hv = H uv^{-1}av$. 
 {}From the preceding lemma, we have $|\tilde{R}_{2n+1} |= (2n+1) 2^n $.

\begin{lemma} \label{qhomlem}
There is a surjective quandle homomorphism $f:\tilde{R}_{2n+1} \rightarrow R_{2n+1}$.  \end{lemma}

The good involution on $\tilde{R}_{2n+1}= (G_{2n+1}, H, a)$ is induced by the map
 $\rho : \tilde {R}_{2n+1} \to \tilde {R}_{2n+1}$ by 
\[
\rho(Hu)=
\left\{
\begin{array}{ll}
H(-1,-2,\ldots, -(n+1),n+2,\ldots, 2n+1)u&\mbox{ if $n$ is an odd number,}\\
H(1,-2,\ldots, -(n+1),n+2,\ldots, 2n+1)u&\mbox{ if $n$ is an even number.}
\end{array}
\right.
\]

Let ${\cal I}$ be the kernel of 
 the map $f: G_{2n+1}\to D_{2(2n+1)}$ 
 in Lemma~\ref{qhomlem}, 
 i.e., $${\cal I}=\{I_{\vec{\epsilon}}=(\epsilon_1(1),\ldots , \epsilon_{2n+1}(2n+1) ) ~|~\epsilon_i =\pm 1, ~\prod _{i=1}^{2n+1}\epsilon _i=1\}.$$
Define the maps $f_a: {\cal I} \to {\cal I}$, $f_b: {\cal I} \to {\cal I}$ and $f_{b^{-1}}: {\cal I} \to {\cal I}$ as follows:
\begin{eqnarray*}
f_a[(\epsilon_1(1),\ldots , \epsilon_{2n+1}(2n+1) )]
&=& (\epsilon_1(1),\epsilon_{2n+1}(2) , \ldots , \epsilon_{2}(2n+1) ), \\ 
f_b[(\epsilon_1(1),\ldots , \epsilon_{2n+1}(2n+1) )]
&=& (\epsilon_{2n+1}(1),\epsilon_1(2),\ldots , \epsilon_{2n}(2n+1) ),\  {\rm and} \\ 
f_{b^{-1}}[(\epsilon_1(1),\ldots , \epsilon_{2n+1}(2n+1) )]
&=& (\epsilon_{2}(1),\ldots ,\epsilon_{2n+1}(2n), \epsilon_{1}(2n+1) ), 
\end{eqnarray*} 
where the order of $(\epsilon_2, \ldots, \epsilon_{2n+1} )$ is reversed for $f_a$, 
and $(\epsilon_1, \ldots, \epsilon_{2n+1} )$ is cyclically permuted for $f_b$ and $f_{b^{-1}}$. 
We can easily see that $f_a$, $f_b$ and $f_{b^{-1}}$ are  
automorphisms of ${\cal I}$, 
$f_a^{-1}=f_a$ and $f_b^{-1}=f_{b^{-1}}$. Moreover, their actions on the diagonal matrices correspond to the dihedral actions of the elements $x$, $y$, and $y^{-1}$. 

We also consider distinguished elements  
\begin{eqnarray*}
I_{+} &=& (1,\ldots,n,-(n+1),n+2,\ldots ,2n, -(2n+1)), \ {\rm and} \\ 
I_{-} &=& (-1, 2,\ldots, n+1, -(n+2), n+3, \ldots ,2n+1), 
\end{eqnarray*} 
which have exactly two minus signs.  
Observe 
that $f_b(I_+)=I_-$. 
Then the following equalities hold: 
$${I_{\vec{\epsilon }}}\; a=af_a(I_{\vec{\epsilon }}), \quad 
{I_{\vec{\epsilon }}}\, b=bf_b(I_{\vec{\epsilon }}), \quad 
 {I_{\vec{\epsilon }}}\, b^{-1}=b^{-1}f_b^{-1}(I_{\vec{\epsilon }}), $$
$$ \ ba=ab^{-1}I_+, \quad
{\rm and}  \quad  \ b^{-1}a=abI_-. $$
Consequently, we have the following relations:
$$b^{\pm j} a =a b^{\mp j} \left[ {\prod}_{k=0}^{j-1} f^{\mp k}_b (I_{\pm}) \right], 
\quad  {\rm and} \quad   
I_{\vec{\epsilon}} \, b^{\pm j} = b^{\pm j} f_b^{\pm j} (I_{\vec{\epsilon}}).$$
If $1 \le j \le n$, then  \vspace{2.5mm} 
${\prod}_{k=0}^{j-1} f^{-k}_b (I_{+})$
is a diagonal matrix 
in $\cal{I}$ that has 
 two blocks of $j$ contiguous $(-)$-signs; the first block ends at  $n+1$, 
and the second block ends at  $2n+1$.  
\vspace{2.5mm} 
In particular, ${\prod}_{k=0}^{n-1} f^{ -k}_b (I_{+})$ 
has exactly one $(+)$-sign at position $(1)$.  \vspace{2.5mm} 
Similarly,
${\prod}_{k=0}^{n-1} f^{ k}_b (I_{-})$ 
has exactly one $(+)$-sign at position $(n+1)$. 

\vspace{0.25cm}  
 
For $i=1, \ldots, 2n+1$, let $I(i)$ denote the diagonal matrix that has 
 exactly one $(+)$-sign at position $(i,i)$ and $(-1)\/$s at the other positions along the diagonal.

\begin{lemma}\label{normal}
The following product formulas hold:
\begin{enumerate}
\setlength{\itemsep}{3pt}
\item  \hfil $(b^i I_{\vec{\epsilon}})(b^j I_{\vec{\delta}}) = b^{i+j} f_b^j(I_{\vec{\epsilon}}) I_{\vec{\delta}},$ \hfill \hfill
\item   \hfil $(b^i I_{\vec{\epsilon}})(a b^j I_{\vec{\delta}}) = ab^{j-i}
\left[ {\prod}_{k=0}^{i-1} f^{j-k}_b (I_{+}) \right] 
 f_b^j(f_a(I_{\vec{\epsilon}})) I_{\vec{\delta}},$ \hfill \hfill
 \item   \hfil $(a b^i I_{\vec{\epsilon}})(a b^j I_{\vec{\delta}}) = a^2 b^{j-i} \left[ {\prod}_{k=0}^{i-1} f^{j-k}_b (I_{+}) \right] f_b^j(f_a(I_{\vec{\epsilon}})) I_{\vec{\delta}}, $ \hfill \hfill  
 \item   \hfil $(a b^i I_{\vec{\epsilon}})( b^j I_{\vec{\delta}}) = a b^{i+j}
f_b^j(I_{\vec{\epsilon}}) I_{\vec{\delta}}.$ \hfill \hfill
\end{enumerate}
\end{lemma}
{\it Proof.} The calculations follow directly from the formulas above. $\Box$

Lemma~\ref{normal}  gives normal forms for the elements of $G_{2n+1}$. We\cite{COS} used these forms to calculate that $\tilde{R}_{2n+1}$ is connected.

A very interesting result therein is that the symmetric homology of $\tilde{R}_3$ contains an infinite cyclic summand.

\section{Applications of Quandles and Quandle Cocycles}

Here, I will give short descriptions of the myriad of results in computations of quandles, their homology and cohomology, and applications of their cocycles and cocycle invariants. 
The contents and organization of this section are based upon a talk that Masahico  Saito gave at Knots in Washington. His original beamer presentation is available at 
\begin{verbatim}
http://shell.cas.usf.edu/~saito/talks/talk09DCquandle.pdf
\end{verbatim}

\noindent
{\bf Quandles with a few elements and their (co)homologies.}
Computations of quandles, biracks, and their (co)-homology have largely been done via computer programs. Some algebraic methods also have been  developed, but the general resolution/derived functor approach is lacking for quandles. 
While analogues between group and quandle cohomology theory are lacking for high degrees, we do know that extensions of quandles by abelian groups are quantified by quandle $2$-cocycles. That is if $\phi: X \times X \rightarrow A$ represents a quandle $2$-cocycle on $X$ that takes values in an Abelian group $A$, then the product $(a,x) \lt ( b, y ) = (a + \phi(x,y), x \lt y)$ defines a quandle structure on $A \times X$. Similar constructions almost always hold in more general contexts\cite{AG,CENS,COS}; see the discussion above.

In the appendix to \cite{CKS:book} a complete list of quandles with up to $6$ elements is given and computations of the quandle homology is provided. However, there are some errors in the homology computations since Uegaki's program did not distinguish between torsion groups of the form $\Z_p \oplus \Z_p$ and $\Z_{p^2}$. These errors have been corrected by Niebrzydowski and Przytycki who used programs that Shumakovich developed originally to compute Khovanov homology.  Przytycki recently sent me the following information:
\begin{itemize}
\item for the quandle $X=QS_4$ that is the last four element quandle listed on p. 172, they computed:
$$H^{\rm Q}_4(X)\cong \Z_4 \oplus \Z_2^4,$$
and
$$H^{\rm Q}_5(X)\cong \Z_4 \oplus \Z_2^5.$$
In addition, they computed $H^{\rm Q}_6(X)\cong \Z_4 \oplus \Z_2^9.$
\item For the quandle $X=QS_6$ (which up to a labeling of elements is the quandle of $4$-cycles in the symmetric group on $4$ elements: 
$$QS_6=\{(1234), (1243), (1324), (1342), (1423), (1432)\}),$$
they compute
$$H^{\rm Q}_3(X) \cong H^{\rm Q}_4(X)=\Z_3\oplus \Z_8.$$
\end{itemize}

Sam Nelson and his students \cite{MN,MNT,CM}
have extended the tables of quandles in a number of directions, and they have provided python code to implement calculations of quandles via the matrix that represents the multiplication table. See for example,
\begin{verbatim}
http://www.esotericka.org/cmc/quandles.html
\end{verbatim}

Gra\~{n}a and his student Vendramin are developing GAP programs to recognize quandles, racks, and to compute their  homology. See
\begin{verbatim}
http://mate.dm.uba.ar/~lvendram/rig/index.html
\end{verbatim}

On the algebraic side, Nelson has an classification of finite Alexander Quandles, and Gra\~{n}a has a classification of indecomposable racks of order $p^2$ \cite{Grana}. The idea  of extensions\cite{CENS} via quandle $2$-cocycles is brought to its full generality\cite{AG}; however, there are still only few examples that are known. The interest in quandles among algebraists has to do with their relationships to Yetter-Drinfeld modules and to solutions to the Yang-Baxter equations \cite{AG,Ei2}. A solution to the set-theoretic YBE is morally equivalent to a biquandle.

Niebrzydowski \cite{Nie} began the problem of classifying Burnside keis.  Joyce's and Matveev's classification of quandles (Theorem~\ref{thm:JM})  as an operation on a coset space has been used\cite{COS} to construct a family of quandles that have interesting properties; see Theorem~\ref{extthm}.

In terms of homology and cohomology there are still very limited results. The appendix to our book\cite{CKS:book} gives calculations for low dimensional quandle homology
for quandles with up to six elements (see corrections above).  We\cite{CEGS} were able to find non-trivial $3$-cocycles, via a computer search, for a quandle in the extended  Andruskiewitsch and Gra\~{n}a theory. Niebrzydowski and Przytycki \cite{NP3} computed up to $H_{12}$ for the dihedral quandle of order $3$. They proposed the ``delayed Fibonacci conjecture'' for $R_p$ namely $H_{\rm Q}^n(R_p) = \Z_p^{f_n}$ where the sequence $f_n$ satisfies the recurrence relation $f_n = f_{n-1} + f_{n-3}$, and $ f(1) = f(2) = 0$, $f(3) = 1$.
 Nosaka has given a proof of this in certain cases. The Japanese summary of his results is found here:
 \begin{verbatim}
www.f.waseda.jp/taniyama/math-of-knots/knots2009proc.pdf
 \end{verbatim}
 
An early result on quandle homology is given by Litherland and Nelson \cite{LN} is that the rank of the integral quandle, rack, and degenerate homology is determined from the number of orbits  of the inner automorphism group of the quandle. Etingof and Gra\~{n}a\cite{EG} prove a similar result and go on to show that the second rack cohomology can be expressed as $H^2(X,A) \cong H^1(G_X, {\mbox{\rm Fun}}(X, A))$ where $G_X$ denotes the enveloping group, and   ${\mbox{\rm Fun}}(X, A)$ denotes the set of functions from the quandle $X$ to the abelian group $A$. Both pairs of authors improved upon our results\cite{CJKS2}.  

Mochizuki \cite{Mochi} computed the second and third homology groups of dihedral quandles of prime order.  Ameur and Saito \cite{AS} found a large number of non-trivial cocycles for Alexander quandles by generalizing the polynomial formulas for Mochizuki cocycles. Some very promising  techniques developed by 
Niebrzydowski and Przytycki  \cite{NP2,NP3} include homology operations to the chain groups that are algebra versions of twist spinning and  partial derivatives defined at the chain level.

Computations of cocycle invariants beyond the original papers \cite{CJKLS:ERA,CJKLS} began with one of our papers\cite{CJKS1}. Another group of collaborators and I also computed invariants associated to the Andruskiewitsch and Gra\~{n}a cocycles\cite{CEGS}. In particular, we showed that many $2$-twist-spun knots were non-invertible using cocycles that we found by computer searches. Nelson and his students have developed more subtle invariants from quandles and biquandles, and they have applied these results to virtual knots and their generalizations (For full bibliographic information on Nelson and his collaborators, please search the ArXiv and MathSciNet). Ameur and Saito \cite{AS} developed used polynomial valued cocycles to study the tangle embedding problem. 
At about the same time, Niebrzydowski  was studying similar problems using cocycle techniques. The link 
\begin{verbatim}
http://shell.cas.usf.edu/quandle/
\end{verbatim} will send you to the database compiled by Masahico's master's degree student Chad Smudde. Included among the associated pages are several Maple programs that are used find and compute cocycles and their invariants. Cocycles are used to show that certain tangles cannot be embedded into specific knots\cite{AERSS}. 

\bigskip

\noindent
{\bf Applications.}
Applications of the quandle $3$-cocycle invariants were developed by various authors in various combinations. Infinite families of non-invertible $2$-knots were found by Asami and Satoh\cite{AsamiSatoh}; Iwakiri\cite{Iwakiri1} used cocycle invariants to determine lower bounds for the number of $1$-handles that are needed to remove all triple points as well as bounds on the triple point number. Computations of the $3$-cocycle invariant of $2$-twist-spun $2$-bridge knots are given\cite{Iwakiri2}. Asami  and Kuga\cite{AK} classify the non-trivial colorings of torus knots by Alexander quandles, and use Satoh's method \cite{Satoh} of twist-spinning to compute  cocycle invariants for selected twist spins of these knots.  Satoh and Shima \cite{SS} used cocycle invariants to show that the $2$-twist-spun trefoil has at least $4$ triple points. These techniques were applied by Hatakanaka \cite{eri} to show that the $2$-twist-spun figure $8$ knot (and the $(2,5)$-torus knot) each have at least $6$-triple points. 

Mohamad and Yashiro compute minimal triple point numbers for twist-spun knots that are colored by the dihedral quandle of order 5. An abstract of that talk is found here:
\begin{verbatim}
http://faculty.ms.u-tokyo.ac.jp/~topology/abstracts.pdf
\end{verbatim}

The $3$-cocycles defined over symmetric quandle homology  (see above) were used by Oshiro and Kamada and Oshiro \cite{O,KO} to determine lower bounds for triple point numbers of linked surfaces in which at least one component is non-orientable. We\cite{COS}  found families of non-orientable connected surfaces in thickened $3$-manifolds whose projections contain as many triple points as you might like (see also the material above). It is quite likely that there are knotted (hence connected) non-orientable surfaces in $4$-space whose projections will have an arbitrarily large number of triple points, and I would guess that symmetric quandle $3$-cocycle can be used to detect them. 

\bigskip

There is also a $3$-cocycle invariant for classical knots. Its definition is similar to the $2$-cocycle invariant, but $2$-dimensional regions are also colored by quandle elements, normals to arcs point from a region colored $a$ to a region colored $a\lt b$ where $b$ is the color on the arc that separates the two regions. The $3$-cocycle invariant  for classical knots can be shown to detect the chirality of the trefoil \cite{FR2}.

In private communications, Shin Satoh has informed us that he can detect chirality of spacial graphs using quandle cocycle invariants. Ishii and Iwakiri \cite{II} have used quandle cocycles to distinguish flowed oriented spacial graphs. In particular, they can distinguish whether or not two knotted spacial graphs represent the same knotted handle body. By using cocycle-like techniques, my collaborators and I were able to show that certain knot isotopies required at least so-many type III Reidemeister moves \cite{CESS}. 

A classical knot is non-trivial if and only if the second quandle cohomology group of the fundamental quandle is non-trivial \cite{Ei}. Eisermann goes on to show more. In general, $H_2(\pi_Q(k)) \cong \Z$ and a choice of generator is tantamount to a choice of orientation class. Thereby, the classification of the knot complement via the quandle is complete. 

Masahico Saito \cite{MS:Fox} showed that Mochizuki's cocycle could be used to give necessary conditions for  $p$ colorable knots (with $p$ a prime number larger than $7$)   to use at least $5$ colors.  Shin Satoh and Masahico Saito \cite{sheetno} used cocycle invariants to give bounds for the minimal number of unbroken sheets in knotted surface diagrams. For example, any diagram of the spun trefoil has at least $4$ sheets.   The three of us used quandle $2$-cocycles to determine that certain surfaces are not ribbon condordant\cite{ribbon}. 

In Japanese, there is an announcement that Kabaya and Inoue relate the quandle cocycle invariants to the Chern-Simon invariants. 
\begin{verbatim}
http://faculty.ms.u-tokyo.ac.jp/~topology/abstracts.pdf
\end{verbatim}
This result is a follow-up to \cite{Inoue} which announces that hyperbolic volume can also be understood as a quandle cocycle invariant. 
Additional topological applications are related to Lefschetz
and elliptic fibrations  via the Dehn quandle \cite{Zablow,KamadaMat}.

 As mentioned above, there is interest among algebraists in quandles and racks because of their relationships to point Hopf algebras, set theoretic solutions to the Yang-Baxter equations, and their roles in Yetter-Drinfeld Modules. In addition, there are connections among categorical groups, Lie $2$-algebras, and other notions of categorification or categorical internalizations. We have been involved in collaborations\cite{CCES1,CCES2} that discuss these connections, and  I anticipate forth-coming work on related topics.

\section{Conclusion}

Racks, quandles, their generalizations and their cohomology theories are finding wide applications in topology and the classifications of algebraic structures. The generalizations are far-reaching because the quandle structure is a natural structure to consider from several points of view. The current paper is a broad-brush painting of some of the current trends and applications that have been found. In as far as I am able, I hope that I pointed the reader to the best, and most up-to-date, bibliographic sources. 

To the authors whose work I neglected to mention, please accept my apologies; the current paper is being prepared in the face of some deadlines, and in the face of other projects of my own  and my collaborators that remain incomplete. I recognize, as should the reader, that the range of results and applications in the age of rapid information dissemination presents an ever-moving target. 
Since there is some lag between time of submission and publication, I will, if possible include summaries of other results if you call my attention to them.

I would like to complete this work with acknowledgements to my many collaborators on these projects, and in particular to my long time collaborator Masahico Saito. When faced with the task of writing this summary, he asked to be excluded only because he too is involved in many other projects. Much of the mathematics that is developed herein is due in large part to his creativity. 
And as I prepared the manuscript, he was very helpful with bibliographic references.

\end{document}